\input ./style/arxiv-general.cfg
\documentclass[MSNbibl,number,citesort,dvips]{arxbj}
\makeatletter
   \@ifpackageloaded{graphicx}{}{\usepackage{graphicx}}
\makeatother

\aid{0}
\volume{21}
\issue{4}
\pubyear{2015}
\firstpage{2073}
\lastpage{2092}
\doi{10.3150/14-BEJ635} 
\docsubty{FLA}

\makeatletter

\newcommand{\rrvert}{\vert}
\newcommand{\llvert}{\vert}
\newcommand{\rrVert}{\Vert}
\newcommand{\llVert}{\Vert}
\newcommand{\lleft}{\left}
\newcommand{\rright}{\right}
\newcommand{\E}{\mathbb{E}}
\newcommand{\p}{\mathbb{P}}
\newcommand{\bfH}{\mathbf{H}}
\newcommand\tbfH{\tilde{\mathbf{H}}}
\newcommand{\bfZ}{\mathbf{Z}}

\newcommand*{\Ind}[1]{\mathbf{1}_{#1}}
\newcommand{\calX}{\mathcal{X}}
\newcommand{\bfX}{\mathbf{X}}
\newcommand{\calF}{\mathcal{F}}
\newcommand{\calZ}{\mathcal{Z}}
\newcommand{\calS}{\mathcal{S}}
\newcommand{\bfy}{\mathbf{y}}
\newcommand{\R}{\mathbb{R}}
\newcommand{\N}{\mathbb{N}}

\newtheorem{lemma}{Lemma}[section]
\newtheorem{theorem}[lemma]{Theorem}
\newtheorem{prop}[lemma]{Proposition}
\newremark{rem}[lemma]{Remark}

\newcommand{\eqref}[1]{(\ref{#1})}

\def\afrac#1#2{#1/(#2)}

\renewcommand{\emptyset}{\varnothing}

\def\afrac#1#2{#1/(#2)}

\def\sklafrac#1#2{(#1/(#2))}

\makeatother

\begin{document}
\begin{frontmatter}

\title{Some remarks on MCMC estimation of spectra of integral operators}
\runtitle{Some remarks on MCMC estimation of spectra of integral operators}

\begin{aug}
\author[A]{\inits{R.}\fnms{Rados{\l}aw}~\snm{Adamczak}\corref{}\thanksref{e1}\ead[label=e1,mark]{R.Adamczak@mimuw.edu.pl}} \and
\author[A]{\inits{W.}\fnms{Witold}~\snm{Bednorz}\thanksref{e2}\ead[label=e2,mark]{W.Bednorz@mimuw.edu.pl}}
\address[A]{Institute of Mathematics, University of Warsaw, ul. Banacha 2, 02-097 Warszawa,
Poland.\\ \printead{e1}; \printead*{e2}}
\end{aug}

%
\received{\smonth{11} \syear{2013}}

%
\begin{abstract}
We prove a law of large numbers for empirical approximations of the
spectrum of a kernel integral operator by the spectrum of random
matrices based on a sample drawn from a Markov chain, which complements
the results by V. Koltchinskii and E. Gin\'{e} for i.i.d. sequences. In
a special case of Mercer's kernels and geometrically ergodic chains, we
also provide exponential inequalities, quantifying the speed of convergence.
\end{abstract}

%
\end{frontmatter}

\section{Introduction}

Let $(\calX,\calF)$ be a measurable space. Consider a probability
measure $\pi$ on $(\calX,\calF)$ and a symmetric measurable kernel
$h\dvtx \calX\times\calX\to\R$, square integrable with respect to
$\pi\otimes\pi$. With $h$ one can associate the kernel linear
operator defined by the formula
%
\begin{eqnarray}
\label{eq:H_def} \bfH f(x) = \int_\calX h(x,y)f(y)\pi(
\mathrm{d}y).
\end{eqnarray}
This is a Hilbert--Schmidt self-adjoint operator on $L_2(\pi)$ and as
such it possesses a real spectrum consisting of a square summable
sequence of eigenvalues. In \cite{GineKoltchinskii}, Koltchinskii and
Gin\'e investigated the problem of approximating the spectrum of $\bfH
$ by the spectra of certain finite dimensional random operators
constructed with the help of the function $h$ and a sequence of i.i.d.
random variables $(X_n)_{n\ge0}$, distributed according to $\pi$.
More precisely, they define a sequence of random matrices
%
\begin{eqnarray}
\label{eq:def_mtH} \tbfH_n = \frac{1}{n} \bigl(h(X_i,X_j)
\bigr)_{0\le i,j\le n-1}
\end{eqnarray}
and
%
\begin{eqnarray}
\label{eq:def_mH} \bfH_n = \frac{1}{n} \bigl( (1-
\delta_{ij})h(X_i,X_j) \bigr)_{0\le i,j \le n-1} =
\tbfH_n - \frac{1}{n}\operatorname{diag} \bigl(
\bigl(h(X_i,X_i) \bigr)_{i=0}^{n-1}
\bigr)
\end{eqnarray}
(above $\delta_{ij}$ is the Kronecker's symbol) and show that with
probability one the spectrum of $\bfH_n$ (completed to an infinite
sequence with zeros) converges in a certain metric to that of $\bfH$.
They also show by simple examples that in general one cannot replace
$\bfH_n$ with $\tbfH_n$. Moreover, under some stronger assumptions,
they provide rates of convergence as well as infinite-dimensional limit
theorems.

Besides intrinsic mathematical interest, the original motivation in
\cite{GineKoltchinskii} came from the limiting theory of
$U$-statistics. A $U$-statistic of degree 2, based on a kernel $h$ and
a sequence $\bfX= (X_n)_{n\ge0}$ is a random variable of the form
%
\begin{eqnarray}
\label{eq:def_Ustat} U_n(h) = U_n(h,\bfX) = \frac{1}{n(n-1)}
\sum_{0\le i\neq j \le n-1} h(X_i,X_j).
\end{eqnarray}
It is well known, that under certain assumptions and proper
normalization, the law of $U_n(h)$ converges to a random variable of
the form $\sum_{i}\lambda_i (g_i^2 - 1)$, where $g_i$'s are i.i.d.
standard Gaussian variables and $\lambda_i$'s are the eigenvalues of
$\bfH$. Thus, the approximate knowledge of the spectrum of $\bfH$
allows for approximate sampling from the limiting spectral distribution
of corresponding $U$-statistics. Since the publication of \cite{GineKoltchinskii}, empirical approximations of spectra found further
applications, for example, in machine learning, especially in the
theory of spectral clustering on manifolds and in the Kernel Principal
Component Analysis (see, e.g., \cite{ShaweTaylorWilliamsChristianiniKandola,SmaleZhou,vonLuxburgBelkinBousquet,RosascoBelkindeVito}).

Although the authors of \cite{GineKoltchinskii} do not develop
specific applications, their results can be interpreted as a Monte
Carlo method for approximating the spectrum of a kernel operator.
However, such an approach would require access to an i.i.d. sample from
the distribution $\pi$, whereas for many situations of interest the
density of the underlying probability measure is known only up to
constants. In such situations, random samples approximating $\pi$ can
be often obtained via Markov Chain Monte Carlo (MCMC) methods, which
rely on simulating a Markov chain with a simple transition function and
invariant measure $\pi$. By the ergodic theorem, after sufficiently
many steps the value of the chain will be distributed approximately as
$\pi$. There are two popular ways of using such samples with
estimators. One of them is to generate sufficiently many independent
samples and to plug them in the estimator. Another one is to use the
estimator directly on the dependent sample coming from the Markov
chain. While the former approach requires analysis of the stability of
the estimated quantity with respect to a small perturbation of the
probability measure, the latter one requires laws of large numbers in
the dependent setting, which would justify using the estimator directly
on the Markov chain.

The objective of this paper is to provide such a law of large numbers,
together with some probability bounds for the problem of approximation
of the spectrum of an integral operator. Our motivation is manifold.
First, we believe that extending the results of Koltchinskii and Gin\'e
to a dependent setting is an interesting probabilistic problem in its
own right. At the same time, it indicates a possibility of having
practical MCMC methods of estimating spectra. Of course, a practical
implementation of this approach would require overcoming additional
obstacles related, for example, to numerical inaccuracy; however, the
law of large numbers and probabilistic bounds provide its theoretical
justification. Additionally, our results
suggest that it should be possible to justify the validity of at least
some of the aforementioned machine learning methods in a dependent
case, which may more accurately model real-life situations.

As a tool, we also develop a law of large numbers for $U$-statistics of
Markov chains started from an arbitrary initial distribution, which
complements results from \cite{ABDGHW,ArconesUstatAbsReg,BertailClemenconUstat,AndrieuJasraDoucetdelMoral}.

The organization of the paper is as follows. First, in Section~\ref{sec:Main_results} we
formulate our results, next in Section~\ref{sec:preliminaries} we present basic notation and preliminary facts
concerning Markov chains (in particular the regeneration method) as
well as tools from linear algebra which will be used in the proofs. In
Section~\ref{sec:Ustat}, we prove the law of large numbers for
$U$-statistics, and in Sections~\ref{sec:Proof_LLN} and \ref
{sec:proof_exp} we provide the proofs of our main results. Finally, in
the last section we discuss the optimality of our assumptions.

\section{Main results}
\label{sec:Main_results}

We will work with a measurable space $(\calX,\calF)$, where $\calF$
is a countably generated $\sigma$-field. Let $\bfX= (X_n)_{n\ge0}$
be a Harris ergodic Markov chain with transition function $P\dvtx
\calX\times\calF\to[0,1]$ and let $\pi$ be its unique invariant
probability measure (we refer to \cite{MT,Numm} for the general theory
of Markov chains on not necessarily countable spaces). We will consider
a symmetric measurable kernel $h\dvtx \calX\times\calX\to\R$ and
the corresponding kernel type operator $\bfH$ given by \eqref
{eq:H_def}. Let $\tbfH_n$ and $\bfH_n$ be random matrices given by
\eqref{eq:def_mtH} and \eqref{eq:def_mH} respectively.

Since the infinite-dimensional operators we will consider will always
be Hilbert--Schmidt, their spectra may be identified with an infinite
sequence $\lambda= (\lambda_n)_{n\ge0} \in\ell_2$, where $\ell_2$
is the Hilbert space of all square summable sequences. There is clearly
some ambiguity here related to the ordering of eigenvalues, but thanks
to the choice of the metric we are about to make, it will not pose a
problem in the sequel, so we may disregard it.

Since we want to approximate the spectrum of $\bfH$ by a spectrum of a
finite-dimensional operator, just as in \cite{GineKoltchinskii} we
will always identify the finite spectrum of the latter with an element
of $\ell_2$, by appending to it an infinite sequence of zeros. We will
denote the spectrum of an operator or a matrix $K$, by $\lambda(K)$.

The metric we will use to compare spectra will be the $\delta_2$
metric defined as
\begin{eqnarray*}
\delta_2(x,y) = \inf_{\sigma\in\mathcal{P}} \Biggl(\sum
_{i=0}^\infty(x_i - y_{\sigma(i)})^2
\Biggr)^{1/2},
\end{eqnarray*}
where $\mathcal{P}$ is the set of all permutations of natural numbers.
It is easy to see that $\delta_2$ is a pseudometric on $\ell_2$.

In what follows, we will always use the notation $\mu f = \int f\,\mathrm{d}\mu$
for a measure $\mu$ and a function~$f$.

Our first result is the following.

%
\begin{theorem}\label{thm:LLN_spectrum} Let $\bfX= (X_n)_{n\ge0}$ be
a Harris ergodic Markov chain on $(\calX,\calF)$ with invariant
probability measure $\pi$ and let $h \dvtx \calX\times\calX\to\R
$ be a symmetric measurable function. Assume that there exists $F\dvtx
\mathcal{X} \to\R$, such that $\pi F^2 < \infty$ and $|h(x,y)| \le
F(x)F(y)$ for all $x,y \in\calX$. Let $\bfH\dvtx  L^2(\pi) \to
L^2(\pi)$ be the linear operator given by \eqref{eq:H_def} and $\tbfH
_n$, $\bfH_n$ be defined by \eqref{eq:def_mtH}, \eqref{eq:def_mH},
respectively. Then for every initial measure $\mu$ of the chain $\bfX
$, with probability one,
\begin{eqnarray*}
\delta_2 \bigl(\lambda(\tbfH_n),\lambda(\bfH) \bigr),
\delta_2 \bigl(\lambda(\bfH _n),\lambda(\bfH) \bigr) \to0.
\end{eqnarray*}
\end{theorem}

Let us now briefly comment on the hypotheses of the above theorem. Our
main assumption is the majorization of the form $|h(x,y)| \le F(x)F(y)$
for some $F\dvtx \calX\to\R$ with $\pi F^2 < \infty$. There are
two main reasons for considering this type of assumptions. The first
one is technical. As shown in \cite{ABDGHW}, the law of large numbers
for $U$-statistics (which we will use in the proofs) of mixing
sequences may fail if one assumes just integrability of the kernel,
which intuitively is related to the fact that the behaviour of the
random variable $h(X_i,X_{i+1})$ may depend on the behaviour of $h$ on
$\pi^{\otimes2}$-negligible sets (since $X_i,X_{i+1}$ are dependent).
As we will see in Section~\ref{sec:examples}, in our setting a similar
phenomenon occurs, in particular the law of large numbers for the
spectra may fail if one assumes only that $\pi^{\otimes2} h^2 <
\infty$. The second reason is the fact that in the theory of Markov
chains, one often proves ergodicity by means of drift conditions and
pointwise assumptions related to the drift functions $V\dvtx \mathcal
{X} \to[0,\infty)$ (see, e.g., \cite{MT,Bax,DFMS,DGM}). The drift
conditions are expressed only in terms of the drift function and the
transition function $P$. While it is not always easy to check
integrability of a general function with respect to the stationary
measure, the drift criteria provide certain integrability for the drift
function. Thus, one can often construct the majorant $F$ in terms of
the function $V$.

Let us also stress that we require that the inequality between $h$ and
$F$ hold pointwise and not just $\pi^{\otimes2}$ a.s. Again, the
reason is related to the dependencies between the variables $X_i$. From
the point of the MCMC applications, it is crucial to allow the Markov
chain to start from arbitrary initial conditions and the distribution
of the chain approaches the stationary measure only in the limit. As a
consequence, it is not enough to assume a $\pi^{\otimes2}$-a.s.
bound. In Section~\ref{sec:examples}, we will illustrate these remarks
with examples.

Finally, let us note that the above theorem provides convergence of
spectra also for the random operator $\tbfH_n$, which as we have
mentioned and as was noted in \cite{GineKoltchinskii} is not the case
in general, even in the i.i.d. setting. To see this, it is enough to
choose a function $h$ vanishing everywhere on $\calX\times\calX$
except for the diagonal, for absolutely continuous $\pi$ and such that
$\int h(x,x)\pi(\mathrm{d}x) = \infty$. The validity of the law of large
numbers for the spectrum of $\tbfH_n$ in our case is of course again a
consequence of our assumptions on $h$ and $F$, which preclude such
counterexamples.

Let us now pass to our second result, which is a tail inequality for
the approximation of spectra. For this, we will work in a more
restrictive, analytic framework, we will also impose stronger
ergodicity assumptions on the chain.

Recall that a Harris ergodic Markov chain with transition function $P$
and invariant measure $\pi$, is geometrically ergodic if there exists
$0< \rho< 1$ such that for every $x \in\calX$ and some constant
$M(x)$, we have for every $n \ge0$,
%
\begin{eqnarray}
\label{eq:geom_erg} \bigl\llVert P^n(x,\cdot)-\pi \bigr\rrVert
_{\mathrm{TV}} \le M(x) \rho^n,
\end{eqnarray}
where $\|\cdot\|_{\mathrm{TV}}$ is the total-variation distance and $P^n$ is
the $n$-step transition function of the chain.

%
\begin{theorem}\label{thm:exponential_ineq} Let $\pi$ be a
probability measure on $(\calX,\calF)$, where $\calX$ is a metric
space and $\calF$ the Borel $\sigma$-field. Let $h\dvtx \calX\times
\calX\to\R$ be a bounded function and $\bfH$ the corresponding
kernel operator defined by \eqref{eq:H_def}.
Assume that there exist continuous functions $\phi_n\dvtx \calX\to
\R$, $n\in I$ (where $I = \{0,\ldots,R\}$ or $I = \N$) which form an
orthonormal system in $L_2(\pi)$ and a sequence of non-negative
numbers $\lambda= (\lambda_n)_{n \in I} \in\ell_2(I)$ such that we
have a point-wise equality
\begin{eqnarray*}
h(x,y) = \sum_{n \in I}\lambda_n
\phi_n(x)\phi_n(y),
\end{eqnarray*}
with the series converging absolutely and almost uniformly on $\calX
\times\calX$.
Assume furthermore that
$\bfX= (X_n)_{n\ge0}$ is a geometrically ergodic Markov chain with
invariant measure $\pi$, started at a point $z$. Then
\begin{eqnarray*}
\p \bigl(\delta_2 \bigl(\lambda(\tbfH_n),\lambda(\bfH)
\bigr) \ge t \bigr) \le2\exp \biggl(-\frac{1}{L}n\min \biggl(
\frac{t^2}{\sup_{x\in\calX}
h(x,x)^2},\frac{t}{\sup_{x\in\calX} h(x,x)} \biggr) \biggr),
\end{eqnarray*}
where the constant $L$ depends only on the transition function $P$ and
the starting point $z$.
\end{theorem}

In the above formulation, we do not specify the dependence of the
constants in the inequality on the parameters of the Markov chain. This
will be done in Section~\ref{sec:proof_exp} via drift conditions.

We state Theorem~\ref{thm:exponential_ineq} for chains started from a
point. In fact, it holds also for chains started from more general
measures $\mu$ satisfying some mild conditions. Since to formulate
this condition we would need to introduce the regeneration technique
for Markov chains, such a formulation is deferred to Remark~\ref
{re:starting_measure} in Section~\ref{sec:proof_exp}.

We remark that the assumptions concerning the function $h$ are
satisfied for continuous positive definite kernels on a large class of
topological spaces. In the case of compact spaces this fact is known as
Mercer's theorem (see, e.g., \cite{MinhNiyogiYao,MR2138444}). Since
there are many generalizations of this result, with subtle differences,
and a discussion of this topic is beyond the scope of this article we
prefer to formulate the theorem in an abstract form.

We remark that similar inequalities in the i.i.d. case were considered,
for example, in \cite{MP_RM_Colt,MP_RMBernoulli} under weaker
assumptions than the boundedness of $h$ (instead some exponential
integrability was assumed). However, those estimates consider a weaker
metric between spectra and, when specialized to the case of bounded
kernels, involve additional logarithmic factors. Thus, Theorem~\ref
{thm:exponential_ineq} (in a version for chains started not necessarily
from a point) improves on their result for bounded kernels even in the
i.i.d. case.

Let us also mention that in our case one can also obtain results for
unbounded kernels, under appropriate drift conditions involving the
function $h$ (using, e.g., results from \cite{AdBed}). However, their
formulation would be much more involved, so we restrict to the special
case of uniformly bounded kernels.

We would like to stress the important role of positive definiteness in
Theorem~\ref{thm:exponential_ineq}. As will be shown in the proof,
thanks to this assumption we can replace the operator $\tbfH_n$ by a
sum of the form $\sum_{i=1}^{n-1} f(X_i)\otimes f(X_i)$ for some
$f\dvtx \calX\to L_2(\pi)$ which is an $L_2(\pi)$-valued additive
functional of the Markov chain $\bfX$ (similar ideas in the i.i.d.
case were used, e.g., in \cite{MP_RM_Colt,MP_RMBernoulli,SmaleZhou}).
This allows to apply the regeneration technique for obtaining
concentration inequalities for additive functionals of Markov chains.

\section{Notation and preliminary facts}
\label{sec:preliminaries}

\subsection{Markov chains}

We will now present basic facts related to the regeneration technique
for Markov chains on general state spaces. This technique was
independently discovered by Nummelin \cite{NummSplit} and Athreya--Ney
\cite{ANsplit} and relies on a decomposition of the trajectory of a
Markov chain into one-dependent paths of random length. Instead of
providing the technical details of the construction, we will just
present its properties, which will be used in the proof. The technical
details can be found in many monographs on Markov chains; we recommend
\cite{MT,Numm,Chen}.

Let thus $(\calX,\calF)$ be a state space, with $\calF$ countably
generated and assume that $P$ is a Markov chain transition function on
$\calX$. Assume also that the corresponding Markov chain $\bfX=
(X_n)_{n\ge0}$ is Harris ergodic. Then there exists a set $C\in\calF
$ with $\pi(C) > 0$, a positive integer $m$, $\delta> 0$ and a
probability measure $\nu$ on $(\calX,\calF)$, such that for all
$x\in C$, $A\in\calF$,
%
\begin{eqnarray}
\label{eq:small_set} P^m(x,A) \ge\delta\nu(A).
\end{eqnarray}

Using the set $C$ for any probability measure $\mu$ one can define two
sequences of random variables $(\tilde{X}_n)_{n\ge0}$, $(Y_n)_{n\ge
0}$ (on some probability space) with the following properties:
\begin{itemize}[(A4)]
\item[(A0)] $(\tilde{X}_n)_{n\ge0}$ is a Markov chain, $\tilde{X}_0
\sim\mu$.
\item[(A1)] $Y_n \in\{0,1\}$.
\item[(A2)] The stopping times $T_0 = \inf\{ k \ge0\dvt  Y_k = 1\}
$, $T_i = \inf\{k> T_{i-1}\dvt  Y_k = 0\}$ are almost surely finite.
Moreover, $T_0,T_1-T_0,T_2-T_1,\ldots$ are independent random
variables, whereas $T_1-T_0,T_2-T_1,\ldots$ are i.i.d. and their
distribution depends only on $P$ (and not on $\mu$). Moreover, $\E
(T_1-T_0) < \infty$.
\item[(A3)] The blocks $Z_i = (\tilde{X}_{m(T_i+1)},\tilde
{X}_{m(T_i+1)+1},\ldots,\tilde{X}_{mT_{i+1} + m -1})$ form a
one-dependent stationary sequence of random variables with values in
$(\calZ,\calS)$, where $\calZ= \bigcup_{k=1}^\infty\mathcal
{X}^k$, $\calS= \sigma(\bigcup_{k=1}^\infty\calF^{\otimes k})$
(i.e., for all $k$, the $\sigma$-fields $\sigma(Z_i\dvt  i<k)$ and
$\sigma(Z_i\dvt  i>k)$ are independent).
\item[(A4)] For any $f \in L_1(\pi)$ and all $k$,
\begin{eqnarray*}
\E\sum_{i=m(T_k+1)}^{mT_{k+1}+m-1} f(\tilde{X}_i)
= m\E(T_1-T_0)\pi f.
\end{eqnarray*}
\end{itemize}

As already mentioned, in the proofs we will use only the above
properties and so we do not present the general construction of the
chain $(\tilde{X}_n)_{n\ge0}$. Let us however briefly describe the
intuition hidden behind it in the special case of $m = 1$. Informally,
if one attempts to generate\vspace*{2pt} the chain then one draws $\tilde{X}_0$
according to the measure $\mu$, and next if at step $n$ one has
$\tilde{X}_n = x$, then for $x \notin C$, the next variable $\tilde
{X}_{n+1}$ is drawn from the distribution $P(x,\cdot)$ and one sets
$Y_n = 0$. If $x \in C$ then one tosses a coin with heads probability
equal to $\delta$. If one gets heads, then $\tilde{X}_{n+1}$ is
generated according to $\nu$ and $Y_n$ is set to one, otherwise $Y_n$
is set to zero and $\tilde{X}_{n+1}$ is generated according to the
probability measure
\begin{eqnarray*}
Q(x,\cdot) = \frac{P(x,\cdot) - \delta\nu(\cdot)}{1-\delta}.
\end{eqnarray*}
It is straightforward but slightly tedious to formalize this intuition
and prove that for Harris ergodic chains it gives properties
(A0)--(A4). For general $m$, one can still repeat this construction for
the $m$-step transition function to define the chain $(\tilde
{X}_{nm})_{n\ge0}$ and then fill in the intermediate variables in such
a way that properties (A0)--(A4) are still satisfied (note that for
$m=1$ the blocks $Z_i$ of property (A3) are in fact independent, which
is not necessarily the case for general $m$). We refer the reader to
\cite{MT,Numm,Chen} for the details.

Since the Markov chain $(X_n)_{n\ge0}$, started from $\mu$ has the
same distribution as $(\tilde{X}_n)_{n\ge0}$ above, to prove a limit
theorem for $(X_n)_{n\ge0}$ it is enough to do it for $(\tilde
{X}_n)_{n\ge0}$ for which one can exploit the additional structure
given by the auxiliary variables $(Y_n)_{n\ge0}$, which often allows
to reduce the proof to the corresponding limit theorem in the
one-dependent or independent case. This strategy has been adopted for
many problems, including the law of large numbers, the central limit
theorem or the law of the iterated logarithm. We again refer to \cite{MT,Numm,Chen} for a detailed exposition.
As a consequence, for the purpose of proving limit theorems, we can
identify the sequences $(X_n)_{n\ge0}$ and $(\tilde{X}_n)_{n\ge0}$.
In what follows, we will adopt this convention (in particular we will
drop the tilde in $\tilde{X}_n$).

In the proofs, we will use the strong law of large numbers for Markov
chains, which can be easily proved using the regeneration method (see
\cite{MT,Numm}).

%
\begin{theorem}\label{thm:LLN_MC} Let $\bfX= (X_n)_{n\ge0}$ be a
Harris ergodic Markov chain on $(\calX,\calF)$, with invariant
probability measure $\pi$ and let $f\dvtx \calX\to\R$ be a $\pi
$-integrable function. Then with probability one, as $n\to\infty$,
\begin{eqnarray*}
\frac{1}{n}\sum_{i=0}^{n-1}
f(X_i) \to\pi f.
\end{eqnarray*}
\end{theorem}

\subsection{Linear algebra}

The main linear-algebraic result we will need is the Hoffman--Wielandt
inequality. To prove the law of large numbers, it will be sufficient to
use its original finite-dimensional version. However, for the
exponential inequality we will use the infinite-dimensional version
proved in \cite{BhatiaElsner}.

%
\begin{theorem}[(Hoffman--Wielandt inequality)]\label{thm:HW}  If $A,B$
are normal Hilbert--Schmidt operators on some Hilbert space, then
\begin{displaymath}
\delta_2 \bigl(\lambda(A),\lambda(B) \bigr) \le\|A-B
\|_{\mathrm{HS}}.
\end{displaymath}
\end{theorem}

\section{Strong law of large numbers for $U$-statistics of~Markov~chains}
\label{sec:Ustat}
Recall the notation \eqref{eq:def_Ustat}. The aim of this section is
to prove the following.

%
\begin{prop}\label{prop:Ustat} Let $\bfX= (X_n)_{n\ge0}$ be a Harris
ergodic Markov chain on $(\calX,\calF)$ with invariant probability
measure $\pi$ and let $h\dvtx \calX\times\calX\to\R$ be a
symmetric measurable function. Assume that there exists a $\pi
$-integrable $F\dvtx \calX\to\R_+$, such that $|h(x,y)| \le
F(x)F(y)$ for all $x,y \in\calX$. Then for every initial probability
$\mu$ of the chain $\bfX$, with probability one,
\begin{eqnarray*}
U_n(h) \to\pi h
\end{eqnarray*}
as $n \to\infty$.
\end{prop}

We remark that in the literature there are several results concerning
laws of large numbers for $U$-statistics under dependence. In \cite{ABDGHW}, such a result is obtained for a class of ergodic stationary
sequences, under assumption of the same nature as ours. However, we
need the above version, since for MCMC applications it is important to
consider Markov chains started from a point (as the very purpose of
MCMC algorithms is to simulate the stationary distribution, which is
not directly accessible). Results of this type have been obtained
recently, for example, in \cite{BertailClemenconUstat,AndrieuJasraDoucetdelMoral}; however, they require higher order
ergodicity of the chain. We would like to add that the results in \cite{BertailClemenconUstat} are not expressed in terms of point-wise bounds
on the kernel $h$ but rather in terms of integrability of certain
functionals on the paths of the Markov chains. Thus, in general they
are not comparable to Proposition~\ref{prop:Ustat}. On the one hand
they may be applicable to kernels which are not bounded by tensor
products, on the other hand the verification of assumptions may be more
difficult.

To prove Proposition~\ref{prop:Ustat}, we will use the following
result which is a simple corollary to Theorem~5.2. in \cite{ABDGHW}
(we remark that this theorem is stated for $\calZ= \R$, but it is
easy to see that its proof works for an arbitrary measurable space).

%
\begin{lemma}\label{le:ABDGHW} Let $\bfZ= (Z_k)_{k\ge0}$ be a
one-dependent stationary sequence of $(\calZ,\calS)$-valued random
variables and let $H\dvtx \calZ^2 \to\R$ be a symmetric measurable
function. Assume that there exists $F\dvtx \calZ\to\R_+$ such that
$|H(x,y)| \le F(x)F(y)$ for all $x,y \in\calZ$ and $\E F(Z_0) <
\infty$. Then with probability one
\begin{eqnarray*}
U_n(H,\bfZ) \to\E H(Z_0,Z_2)
\end{eqnarray*}
as $n \to\infty$.\vadjust{\goodbreak}
\end{lemma}

\begin{pf*}{Proof of Proposition~\ref{prop:Ustat}}
Define $N_n = \sup\{ k\dvt  mT_k +m-1\le n-1\}$ (with the convention
that $\sup\emptyset= 0$). By the law of large numbers and property
(A2), we have as $n \to\infty$,
%
\begin{eqnarray}
\label{eq:ratio_of_ns} \frac{n}{N_n} \to m \E(T_1-T_0),
\qquad \p_\mu\mbox{-a.s. }
\end{eqnarray}
Recall the space $\calZ$ defined in property (A3). In what follows, we
will use the following convention regarding its elements: for $\mathbf
{x}= (x_1,\ldots,x_k) \in\calZ$ we set $|\mathbf{x}| = k$. Let $H
\dvtx \calZ\times\calZ\to\R$ be the kernel defined by
\begin{eqnarray*}
H(\mathbf{x},\bfy) = \sum_{i=1}^{|\mathbf{x}|} \sum
_{j=1}^{|\bfy
|} h(x_i,y_j).
\end{eqnarray*}
Note that for $\tilde{F} \dvtx \calZ\to\R_+$, given by
\begin{eqnarray*}
\tilde{F}(\mathbf{x}) = \sum_{i=1}^{|\mathbf{x}|}
F(x_i),
\end{eqnarray*}
we have $|H(\mathbf{x},\bfy)| \le\tilde{F}(\mathbf{x})\tilde
{F}(\bfy)$.
Moreover, by property (A4) we have
$\E\tilde{F}(Z_0) < \infty$.

By properties (A3), (A4) and the Fubini theorem, we also get
\begin{eqnarray*}
\E H(Z_k,Z_l) = \bigl(m \E(T_1-T_0)
\bigr)^2 \pi^{\otimes2} h
\end{eqnarray*}
if $|k-l| \ge2$.

Thus, by Lemma~\ref{le:ABDGHW}, \eqref{eq:ratio_of_ns} and the above
equality, we get
%
\begin{eqnarray}
\label{eq:convergence_main} \frac{1}{n(n-1)}\sum_{0\le i\neq j \le N_n-1}
H(Z_i,Z_j) \to\pi ^{\otimes2} h.
\end{eqnarray}

Define also $\tilde{H} \dvtx \calZ\times\calZ\to\R$ as $\tilde
{H}(\mathbf{x},\bfy) = \sum_{i=1}^{|\mathbf{x}|}\sum_{j=1}^{|\bfy
|} |h(x_i,y_j)|$.

In view of \eqref{eq:convergence_main}, to prove the proposition it
remains to show that with probability one the sequences
\begin{eqnarray*}
I_n &=& \frac{1}{n(n-1)}\sum_{i=0}^{mT_0+m-1}
\sum_{j=m(T_0+1)}^{n-1} \bigl|h(X_i,X_j)\bigr|,
\\
\mathit{II}_n & = &\frac{1}{n(n-1)} \sum
_{i=0}^{N_n} \tilde{H}(Z_i,Z_i),
\\
\mathit{III}_n & =& \frac{1}{n(n-1)}\sum
_{i=0}^{N_n-1} \tilde{H}(Z_{N_n},Z_i),
\end{eqnarray*}
converge a.s. to $0$ as $n \to\infty$.

Note that
\begin{eqnarray*}
I_n \le\sum_{i=0}^{mT_0+m-1}F(X_i)
\frac{1}{n(n-1)} \sum_{j=0}^{n-1}
F(X_j) \to0\qquad \mbox{a.s.}
\end{eqnarray*}
since by Theorem~\ref{thm:LLN_MC}, $n^{-1}\sum_{j=0}^{n-1} F(X_j) \to
\pi F$ a.s.

As for $\mathit{II}_n$, we have
\begin{eqnarray*}
\E\tilde{H}(Z_i,Z_i)^{1/2} \le\E\sum
_{i= m(T_i+1)}^{mT_{i+1} + m
-1} F(X_i) = m
\E(T_1-T_0) \pi F < \infty,
\end{eqnarray*}
where we again used (A4).

Thus, using (A3) and \eqref{eq:ratio_of_ns} we get by the
Marcinkiewicz law of large numbers that $\mathit{II}_n \to0$ a.s.

To prove that $\mathit{III}_n\to0$ a.s., note that
\begin{eqnarray*}
\mathit{III}_n = \frac{N_n(N_n+1)}{2n(n-1)}U_{N_n+1}(\tilde{H},\bfZ) -
\frac
{N_n(N_n-1)}{2n(n-1)}U_{N_n}(\tilde{H},\bfZ).
\end{eqnarray*}
By Lemma~\ref{le:ABDGHW} and \eqref{eq:ratio_of_ns}, both terms on
the right-hand side above converge a.s. to
\begin{eqnarray*}
2^{-1} \bigl(m\E(T_1-T_0)
\bigr)^{-2} \E\sum_{i=m(T_0+1)}^{mT_1+m-1} \sum
_{j=m(T_2+1)}^{mT_3+m-1} \bigl|h(X_i,X_j)\bigr|
\le 2^{-1} (\pi F)^2 < \infty,
\end{eqnarray*}
where in the first inequality we used the assumption on $h$ and $F$
together with (A3), (A4) and the Fubini theorem.

This shows that indeed $\mathit{III}_n\to0$ a.s. and ends the proof of
Proposition~\ref{prop:Ustat}.
\end{pf*}

\section{Proof of Theorem \texorpdfstring{\protect\ref{thm:LLN_spectrum}}{2.1}}
\label{sec:Proof_LLN}

To prove Theorem~\ref{thm:LLN_spectrum}, we will need one more simple
result, namely a Marcinkiewicz--Zygmund-type law of large numbers for
Markov chains. Its proof is a standard application of the regeneration
technique. Since we have not been able to find it in the literature, we
provide it for completeness.

%
\begin{lemma}[(Marcinkiewicz--Zygmund LLN for Markov chains)]\label{le:MZ_LLN}
Let $\bfX= (X_n)_{n\ge0}$ be a Harris ergodic Markov chain on $(\calX
,\calF)$ and let $f \dvtx \calX\to\R$ be a measurable function.
Consider $p \in(0,1)$ and assume that $\pi|f|^{p} < \infty$. Then
for any initial measure of the chain, with probability one
\begin{eqnarray*}
\frac{1}{n^{1/p}}\sum_{i=0}^{n-1}
f(X_i) \to0.
\end{eqnarray*}
\end{lemma}

\begin{pf}
As in the proof of Proposition~\ref{prop:Ustat}, define $N = N_n =
\sup\{ k\dvt  mT_k +m-1\le n-1\}$ and recall \eqref{eq:ratio_of_ns}.
Define a function $F\dvtx \calZ\to\R$ (where $\calZ$ is defined in
property (A3) of Section~\ref{sec:preliminaries}) with the formula
\begin{eqnarray*}
F(\mathbf{x}) = \sum_{i=1}^{|\mathbf{x}|}
\bigl|f(x_i)\bigr|.
\end{eqnarray*}
Then, by concavity of the function $t\mapsto|t|^p$ and property (A4)
we get for $i \ge1$,
%
\begin{eqnarray}
\label{eq:integrability_MZ} \E F(Z_i)^p \le\E\sum
_{i=m(T_i+1)}^{mT_{i+1}+m-1} \bigl|f(X_i)\bigr|^p = m
\bigl(\E (T_1-T_0) \bigr)\pi|f|^p < \infty.
\end{eqnarray}
Now
\begin{eqnarray*}
\frac{1}{n^{1/p}} \Biggl\llvert \sum_{i=0}^{n-1}
f(X_i) \Biggr\rrvert \le\frac
{1}{n^{1/p}}\sum
_{i=0}^{mT_0 + m -1} \bigl|f(X_i)\bigr| + \frac{1}{n^{1/p}}
\sum_{i=0}^{N_n} F(Z_i).
\end{eqnarray*}
The first term on the right-hand side above converges a.s. to zero as
$n\to\infty$. Moreover, since $Z_i$ form a stationary one-dependent
sequence by \eqref{eq:ratio_of_ns}, \eqref{eq:integrability_MZ} and
the classical Marcinkiewicz--Zygmund LLN, the second term also
converges a.s. to zero, which ends
the proof of the lemma.
\end{pf}

The proof of Theorem~\ref{thm:LLN_spectrum} will mimic closely the
corresponding proof by Koltchinskii and Gin\'e, in fact one could keep
the linear-algebraic part exactly the same, while replacing just the
probabilistic ingredients (using Proposition~\ref{prop:Ustat} and
Lemma~\ref{le:MZ_LLN}). However, we will slightly change the
exposition with respect to \cite{GineKoltchinskii}, which will allow
to shorten the proof a little bit.

\begin{pf*}{Proof of Theorem~\ref{thm:LLN_spectrum}}
Let us first notice that thanks to the Hoffman--Wielandt inequality and
the assumption on $h$, we have
\begin{displaymath}
\delta_2 \bigl(\lambda(\bfH_n),\lambda(
\tbfH_n) \bigr)^2 \le\llVert \bfH_n - \tbfH
_n\rrVert _{\mathrm{HS}}^2 \le\frac{1}{n^2}\sum
_{i=1}^n F(X_i)^4.
\end{displaymath}
Since $\pi F^2 < \infty$, by Lemma~\ref{le:MZ_LLN} applied with $p =
1/2$, the right-hand side above converges a.s. to zero. Thus, it is
enough to prove the theorem for the matrix $\tbfH_n$.

Since $\pi^{\otimes2} h^2 \le(\pi F^2)^2 < \infty$, $\bfH$ is a
Hilbert--Schmidt operator and so, by the spectral theorem, there exists
an orthonormal system $(\phi_i)_{i \in I}$ in $L_2(\pi)$ (where $I =
\{0,\ldots,R\}$ for some $R \in\N$ or $I =\N$) and a square
summable sequence $(\lambda_i)_{i\in I}$ with non-increasing absolute
values such that
%
\begin{eqnarray}
\label{eq:spectral_representation} h(x,y) = \sum_{i\in I}
\lambda_i \phi_i(x)\phi_i(y),
\end{eqnarray}
where the equality holds in the $L_2(\pi^{\otimes2})$ sense.

As in \cite{GineKoltchinskii} assume first that $h(x,y) = \sum_{i=0}^R \lambda_i \phi_i(x)\phi_i(y)$ and the equality holds
pointwise. Define for $n \ge0$, the sequence of vectors in $\R^n$,
\begin{eqnarray*}
\Phi_i^n = \biggl(\frac{\phi_i(X_0)}{\sqrt{n}},\ldots,
\frac{\phi
_i(X_{n-1})}{\sqrt{n}} \biggr),\qquad 0\le i \le R
\end{eqnarray*}
and note that for $u \in\R^n$,
\begin{eqnarray*}
\tbfH_n u = \sum_{i=0}^R
\lambda_i \bigl\langle\Phi_i^n ,u \bigr
\rangle\Phi_i^n,
\end{eqnarray*}
where $\langle\cdot,\cdot\rangle$ denotes the standard inner
product in $\R^n$.

Now consider the space $\R^{R+1}$ with the standard basis $e_0,\ldots
,e_R$ and let $A_n \dvtx \R^{R+1} \to\R^n$ be the operator given by
$A_n e_i = \Phi_i^n$, $i = 0,\ldots,R$. Define also an operator $K$
on $\R^{R+1}$ as
\begin{eqnarray*}
Ku = \sum_{i=0}^R \lambda_i
\langle e_i,u\rangle e_i.
\end{eqnarray*}
Then, as one can easily check,
\begin{eqnarray*}
\tbfH_n = A_n K A_n^T
\end{eqnarray*}
and since for any two operators $K_1\dvtx \R^a \to\R^b$ and $K_2
\dvtx \R^b \to\R^a$, the (algebraic) spectra of $K_1K_2$ and
$K_2K_1$ are the same (recall our convention of completing the spectra
with zeros to an infinite sequence), we get $\lambda(\tbfH_n) =
\lambda(K A_n^T A_n)$. Together with the obvious equality $\lambda(K)
= \lambda(\bfH)$, this gives
%
\begin{eqnarray}
\label{eq:distance_equality} \delta_2 \bigl(\lambda(\tbfH_n),\lambda(
\bfH) \bigr)= \delta_2 \bigl(\lambda \bigl(KA_n^TA_n
\bigr),\lambda(K) \bigr).
\end{eqnarray}
But for each $i,j = 0,\ldots,R$ we have $\langle A_n^T A_n e_i,
e_j\rangle= \langle A_n e_i,A_n e_j\rangle= \frac{1}{n}\sum_{k=0}^{n-1} \phi_i(X_k)\phi_j(X_k)$. Thus, by Theorem~\ref
{thm:LLN_MC} with probability one, $\langle A_n^T A_n e_i, e_j\rangle
\to\pi\phi_i\phi_j = \delta_{ij}$ and thus\vspace*{1pt} $KA_n^T A_n \to K$,
which implies that the right-hand side of \eqref{eq:distance_equality}
converges to zero a.s. (note that $KA_n^TA_n$ in general is not a
normal matrix, so we cannot use the Hoffman--Wielandt inequality, but
we are working now in a fixed dimension $R+1$ and so we can simply use
the fact that the eigenvalues are continuous functions of the matrix
entries; see, e.g., Appendix D in \cite{HJ}). This proves the theorem
in the special case of finite dimensional kernels.

Consider now an arbitrary kernel $h$, satisfying \eqref
{eq:spectral_representation}. Fix $\varepsilon> 0$. Since $\sum_{i\in
I} \lambda_i^2 < \infty$, there exists $R$ such that $\sum_{i\in I,
i> R} \lambda_i^2 < \varepsilon$. Set $h_R(x,y) = \sum_{i=0}^R
\lambda_i \phi_i(x)\phi_i(y)$ (by which we mean that the equality
holds pointwise, for some particular fixed choice of representatives
from the equivalence class of $\phi_i$ in $L^2(\pi)$). Let $\bfH^R$
be the kernel operator corresponding to $h_R$ and $\tbfH^R_{n} =
(h_R(X_i,X_j))_{0\le i,j\le n-1}$. Define moreover $\tilde{h}_R = h - h_R$.
We have
%
\begin{eqnarray}
\label{eq:approx_oper} \delta_2 \bigl(\lambda(\bfH),\lambda \bigl(
\bfH^R \bigr) \bigr)^2 = \sum
_{i\in I, i >
R}^\infty \lambda_i^2 <
\varepsilon.
\end{eqnarray}
Define the function $F_1 = F + \sum_{i=0}^R \sqrt{|\lambda_i|} |\phi
_i|$ (again we interpret this equality in the pointwise sense) and note
that $F_1 \in L_2(\pi)$. Moreover, for all $x,y\in\calX$,
$|h_R(x,y)|, |\tilde{h}_R(x,y)| \le F_1(x)F_1(y)$.
Thus, by the first part of the proof, we get
%
\begin{eqnarray}
\label{eq:approx_partial} \delta_2 \bigl(\lambda \bigl(\bfH^R
\bigr), \lambda \bigl(\tbfH^R_n \bigr) \bigr) \to0 \qquad
\mbox{a.s.},
\end{eqnarray}
while by Proposition~\ref{prop:Ustat} and Lemma~\ref{le:MZ_LLN} we
obtain that with probability one,
\begin{eqnarray*}
\lim_{n \to\infty} \bigl\llVert \tbfH_n -
\tbfH^R_{n} \bigr\rrVert _{\mathrm{HS}}^2 &=&
\lim_{n\to
\infty} \frac{1}{n^2}\sum_{ i,j =0}^{n-1}
\tilde{h}_R(X_i,X_j)^2
\\
&\le&\lim_{n\to\infty} U_n(h,\bfX) + \lim
_{n\to\infty} \frac
{1}{n^2}\sum_{i=0}^{n-1}
F_1(X_i)^4= \pi^{\otimes2}
\tilde{h}_R^2 = \sum_{i\in I, i> R}^\infty
\lambda_i^2 < \varepsilon.
\end{eqnarray*}
Thus, by the Hoffman--Wielandt inequality,
\begin{eqnarray*}
\limsup_{n \to\infty} \delta_2 \bigl(\lambda(
\tbfH_{n}),\lambda \bigl(\tbfH ^R_{n} \bigr)
\bigr) \le\varepsilon^{1/2} \qquad \mbox{a.s.}
\end{eqnarray*}
In combination with \eqref{eq:approx_oper} and \eqref
{eq:approx_partial}, this implies that for every $\varepsilon> 0$,
\begin{eqnarray*}
\limsup_{n \to\infty} \delta_2 \bigl(\lambda(\bfH),
\lambda( \tbfH _{n}) \bigr)\le2\varepsilon^{1/2} \qquad
\mbox{a.s.}
\end{eqnarray*}
and in consequence
\begin{eqnarray*}
\delta_2 \bigl(\lambda(\bfH),\lambda(\tbfH_{n}) \bigr)
\to0\qquad \mbox{a.s.}
\end{eqnarray*}
\upqed
\end{pf*}

\section{Proof of Theorem \texorpdfstring{\protect\ref{thm:exponential_ineq}}{2.2}}\vspace*{-9pt}
\label{sec:proof_exp}

\begin{pf*}{Proof of Theorem~\ref{thm:exponential_ineq}}
In what follows by $\langle\cdot,\cdot\rangle$ we will denote both
the inner product in $L_2(\pi)$ and in finite-dimensional spaces,
since the precise meaning will always be clear from the context, this
should not lead to ambiguity. The letters $C,c$ will denote absolute
positive constants, whose values may differ between occurrences.

Define $f\dvtx \calX\to L_2(\pi)$ with the formula
\begin{eqnarray*}
f(x) = \sum_{i\in I} \sqrt{\lambda_i}
\phi_i(x)\phi_i(\cdot).
\end{eqnarray*}
Note that $\sum_{i\in I} (\sqrt{\lambda_i}\phi_i(x))^2 = h(x,x) <
\infty$ and that $\phi_i$ form an orthonormal system in $L_2(\pi)$,
so the above series indeed converges in $L_2(\pi)$. Consider now a
random operator on $L_2(\pi)$ given by
\begin{eqnarray*}
K_n = \frac{1}{n}\sum_{i=0}^{n-1}
f(X_i)\otimes f(X_i),
\end{eqnarray*}
that is, for all $u \in L_2(\pi)$
\begin{eqnarray*}
K_n u = \frac{1}{n}\sum_{i=0}^{n-1}
\bigl\langle f(X_i), u \bigr\rangle f(X_i).
\end{eqnarray*}

Note that $K_n$ can be written as $A_nA_n^T$, where $A_n\dvtx \R^{n}
\to L_2(\pi)$ is defined by $A_ne_i = n^{-1/2}f(X_i)$ ($e_0,\ldots
,e_{n-1}$ being the standard basis in $\R^n$). Thus $\lambda(K_n) =
\lambda(A_n^TA_n)$ (recall that we append spectra of finite
dimensional operators with infinite sequences of zeros). But
\begin{eqnarray*}
\bigl\langle A_n^TA_ne_i,e_j
\bigr\rangle= \langle A_ne_i,A_ne_j
\rangle= \frac
{1}{n} \bigl\langle f(X_i),f(X_j)
\bigr\rangle= \frac{1}{n}\sum_{k\in I}\lambda
_k\phi_k(X_i)\phi_k(X_j)
= \frac{1}{n}h(X_i,X_j),
\end{eqnarray*}
so $A_n^TA_n = \tbfH_n$. Thus, our goal will be to bound the distance
between the spectrum of $K_n$ and the sequence $\lambda$.

The random operator $K_n$ is a sum of independent random rank one
operators, moreover, using the fact that $\phi_i$ form an orthonormal
system in $L_2(\pi)$ one easily checks that
%
\begin{eqnarray}
\label{eq:L_infty_bound} \bigl\llVert f(x)\otimes f(x) \bigr\rrVert _{\mathrm{HS}} =
h(x,x)
\end{eqnarray}
and
%
\begin{eqnarray}
\label{eq:mean_operator} \E_\pi f(X_i)\otimes f(X_i)
= \bfH
\end{eqnarray}
(where the expectation on the left-hand side is the Bochner integral in
the Hilbert space of Hilbert--Schmidt operators).

Thus, we can apply to $K_n$ classical results concerning concentration
for sums of independent Banach space valued random variables [after
passing to the block decomposition given by (A3)]. The inequality we
will use is a version of Bernstein's $\psi_1$ inequality. To formulate
it, let us first recall the definition of the Orlicz $\psi_1$ norm.
For a Banach space valued random variable $X$, we define
\begin{eqnarray*}
\llVert X\rrVert _{\psi_1} = \inf \bigl\{\rho>0\dvt \E\exp \bigl(\llVert
X \rrVert /\rho \bigr) \le2 \bigr\}.
\end{eqnarray*}
By exponential Chebyshev's inequality, we have
%
\begin{eqnarray}
\label{eq:Chebyshev} \p\bigl(|X| \ge t\bigr) \le2\exp \bigl(-t/\llVert X\rrVert _{\psi_1}
\bigr)
\end{eqnarray}
for $t > 0$.

The following inequality is a simple corollary to Theorem~1.4. in \cite{VershyninFrame}.

%
\begin{lemma}\label{le:Bernstein_psi_one} Let $U,U_i$, $i = 1,\ldots, n$, be i.i.d. mean zero
random variables with values in a Banach space $(B,\llVert \cdot\rrVert )$. Assume
that $\llVert U\rrVert _{\psi_1} < \infty$. Then for all $t > 0$,
\begin{eqnarray*}
\p \Biggl( \Biggl\llvert \Biggl\llVert \sum_{i=1}^n
U_i \Biggr\rrVert - \E \Biggl\llVert \sum
_{i=1}^nU_i \Biggr\rrVert \Biggr\rrvert
\ge t \Biggr)\le2\exp \biggl(-c\min \biggl(\frac{t^2}{n\llVert U\rrVert _{\psi
_1}^2}, \frac{t}{\llVert U\rrVert _{\psi_1}}
\biggr) \biggr),
\end{eqnarray*}
where $c > 0$ is a universal constant.
\end{lemma}

It is well known (see, e.g., \cite{MT}, Chapters~15, 16, or \cite{RobRos,Bax})
that for uniformly ergodic Markov chains we have $\llVert T_1-T_0\rrVert _{\psi_1}
< \infty$ and if the chain is started from a point, then also $\llVert T_0\rrVert
_{\psi_1}<\infty$, which allows for the use of the above inequality
in our setting.

Let us now define $g(x) = f(x)\otimes f(x)$, $U_i = \sum_{i =
m(T_i+1)}^{mT_{i+1}+m-1} (g(X_i)-\pi g)$ and recall the definition $N =
N_n = \sup\{ k\dvt  mT_k +m-1\le n-1\}$. Using properties (A0)--(A4)
and \eqref{eq:mean_operator} we get that $\E U_i = 0$ and
\begin{eqnarray*}
\Biggl\llVert \frac{1}{n}\sum_{i=0}^{n-1}
g(X_i) - \bfH \Biggr\rrVert _{\mathrm{HS}} &=& \frac{1}{n}
\Biggl\llVert \sum_{i=0}^{n-1}
\bigl(g(X_i) - \pi g \bigr) \Biggr\rrVert _{\mathrm{HS}}
\\
&\le&\frac{1}{n} \Biggl\llVert \sum_{i=0}^{(mT_0 + m - 1)\wedge(n-1)}
\bigl(g(X_i) - \pi g \bigr) \Biggr\rrVert _{\mathrm{HS}} +
\frac{1}{n} \Biggl\llVert \sum_{i=0}^{N-1}
U_i \Biggr\rrVert _{\mathrm{HS}}
\\
&&{} + \frac
{1}{n} \Biggl\llVert \sum_{i=m(T_N+1)}^{n-1}
g(X_i) - \pi g \Biggr\rrVert _{\mathrm{HS}}
\nonumber
\\
&\le&\frac{2m}{n}(T_0 + 1)\llVert g\rrVert _\infty+
\frac{1}{n} \Biggl\llVert \sum_{i=0}^{N-1}
U_i \Biggr\rrVert _{\mathrm{HS}} + \frac{2}{n} \bigl(n-
m(T_N + 1) \bigr)_+\llVert g\rrVert _\infty ,
\nonumber
\end{eqnarray*}
where $\llVert g\rrVert _\infty= \sup_{x\in\calX}\llVert g(x)\rrVert _{\mathrm{HS}}$.

Therefore,
%
\begin{eqnarray}
\label{eq:decomposition} &&\p \Biggl( \Biggl\llVert \frac{1}{n}\sum
_{i=0}^{n-1} g(X_i) - \bfH \Biggr\rrVert
_{\mathrm{HS}} \ge t \Biggr)
\nonumber
\\
&&\quad \le\p \bigl(2m(T_0+1)\llVert g\rrVert _\infty\ge
tn/3 \bigr) + \p \Biggl( \Biggl\llVert \sum_{i=0}^{N-1}
U_i \Biggr\rrVert _{\mathrm{HS}} \ge tn/3 \Biggr)
\\
&&\qquad {}+ \p \bigl( 2 \bigl(n- m(T_N + 1) \bigr)\llVert g\rrVert
_\infty\ge tn/3 \bigr).
\nonumber
\end{eqnarray}

By \eqref{eq:Chebyshev},
%
\begin{eqnarray}
\label{eq:first_part} \p \bigl(2m(T_0+1)\llVert g\rrVert _\infty
\ge tn/3 \bigr) \le2\exp \biggl(-\frac{nt}{6 m\llVert
T_0+1\rrVert _{\psi_1}\llVert g\rrVert _\infty} \biggr).
\end{eqnarray}

Moreover, by Lemma~3 in \cite{AdMarkovTail} we have for all $t > 0$,
\begin{eqnarray*}
\p \bigl(2 \bigl(n- m(T_N + 1) \bigr)\ge t \bigr) \le2\exp
\biggl(-c \frac{t}{m\tau\log\tau
} \biggr),
\end{eqnarray*}
where $\tau= \max\{\llVert T_0+1\rrVert _{\psi_1},\llVert T_1-T_0\rrVert _{\psi_1}\}$ (we
remark that the notation and the definition of splitting times in \cite{AdMarkovTail} are slightly different than ours, in particular the
Markov chain there is indexed by $\N\setminus\{0\}$ and not by $\N$,
however it is easy to see that the simple proof of Lemma~3 can be
carried over to our setting).
Thus,
%
\begin{eqnarray}
\label{eq:third_term} \p \bigl( 2 \bigl(n- m(T_N + 1) \bigr)\llVert g
\rrVert _\infty\ge tn/3 \bigr) \le2\exp \biggl(-c\frac
{nt}{m\llVert g\rrVert _\infty\tau\log\tau}
\biggr).
\end{eqnarray}

To handle the middle term in the decomposition \eqref
{eq:decomposition}, we will apply Lemma~\ref{le:Bernstein_psi_one} to
the random variables $U_i$. Since for $m > 1$, these variables are only
one-dependent; moreover, the number of full blocks $Z_i$ in the
sequence $X_0,\ldots,X_{n-1}$ is random, there are two technical
steps, which have to be carried out first, namely we have to split the
sum $\sum_{i=0}^{N-1} U_i $ into odd and even terms and use a L\'evy
type inequality to handle the random number of summands. We have
\begin{eqnarray*}
&&\p \Biggl( \Biggl\llVert \sum_{i=0}^{N-1}
U_i \Biggr\rrVert _{\mathrm{HS}} \ge t/3 \Biggr)
\\
&&\quad \le\p \biggl( \biggl\llVert \sum_{0\le i \le N-1, 2|i}U_i
\biggr\rrVert _{\mathrm{HS}} \ge t/6 \biggr) + \p \biggl( \biggl\llVert \sum
_{0\le i \le N-1, \neg2| i}U_i \biggr\rrVert _{\mathrm{HS}}
\ge t/6 \biggr)
\\
&&\quad \le C\p \Biggl( \Biggl\llVert \sum_{i=0}^{\lfloor n/m\rfloor-1}
\tilde{U}_i \Biggr\rrVert _{\mathrm{HS}}\ge t/C \Biggr),
\end{eqnarray*}
where $C$ is a universal constant and $\tilde{U}_i$ is a sequence of
independent random variables, distributed as $U_0$. In the last
inequality, we used the fact that $N m \le n$ and a L\'evy-type
inequality for i.i.d. Banach-space valued random variables due to
Montgomery--Smith \cite{MontgomerySmithLevy}, which asserts that for a
sequence $W_i$ of i.i.d. Banach space-valued variables
\begin{eqnarray*}
\p \Biggl(\max_{k\le n} \Biggl\llVert \sum
_{i=1}^k W_i \Biggr\rrVert \ge t
\Biggr) \le C\p \Biggl( \Biggl\llVert \sum_{i=1}^n
W_i \Biggr\rrVert \ge t/C \Biggr).
\end{eqnarray*}

Now, we have
\begin{eqnarray*}
\llVert U_i\rrVert _{\psi_1} \le 2 m\llVert g\rrVert
_\infty \llVert T_{1} -T_0\rrVert
_{\psi_1}
\end{eqnarray*}
and so Lemma~\ref{le:Bernstein_psi_one} gives
\begin{eqnarray*}
&&\p \Biggl( \Biggl\llVert \sum_{i=0}^{N-1}
U_i \Biggr\rrVert _{\mathrm{HS}} \ge C\E \Biggl\llVert \sum
_{i=0}^{\lfloor
n/m\rfloor-1} \tilde{U}_i \Biggr
\rrVert _{\mathrm{HS}} + s/3 \Biggr)
\\
&&\quad \le2\exp \biggl(-c\min \biggl(\frac{s^2}{nm\llVert g\rrVert _\infty^2 \llVert T_1-T_0\rrVert
_{\psi_1}^2},\frac{s}{m\llVert g\rrVert _\infty\llVert T_1-T_0\rrVert _{\psi_1}}
\biggr) \biggr).
\end{eqnarray*}

Using the above bound together with \eqref{eq:decomposition}, \eqref
{eq:first_part} and \eqref{eq:third_term} we arrive (after adjusting
the constants) at
\begin{eqnarray*}
&&\p \Biggl( \Biggl\llVert \frac{1}{n}\sum_{i=0}^{n-1}
g(X_i) - \bfH \Biggr\rrVert _{\mathrm{HS}} \ge Cn^{-1}\E
\Biggl\llVert \sum_{i=0}^{\lfloor n/m\rfloor-1}
\tilde{U}_i \Biggr\rrVert _{\mathrm{HS}} + t \Biggr)
\\
&&\quad \le 2\exp \biggl(-cn\min \biggl(\frac{t^2}{m\llVert g\rrVert _\infty^2 \tau
^2},\frac{t}{m\llVert g\rrVert _\infty\tau\log\tau}
\biggr) \biggr).
\end{eqnarray*}

Using the fact that the norm $\llVert \cdot\rrVert _{\mathrm{HS}}$ is Hilbertian and $\E
\tilde{U}_i = 0$, we obtain
\begin{eqnarray*}
\E \Biggl\llVert \sum_{i=0}^{\lfloor n/m\rfloor-1}
\tilde{U}_i \Biggr\rrVert _{\mathrm{HS}}^2 = \sum
_{i=0}^{\lfloor n/m\rfloor-1} \E\llVert \tilde{U}_i
\rrVert _{\mathrm{HS}}^2 \le\frac
{4n}{m} m^2
\E(T_1 - T_0)^2 \llVert g \rrVert
_\infty^2 \le Cnm\tau^2\llVert g\rrVert
_\infty^2,
\end{eqnarray*}
which combined with the previous inequality gives
\begin{eqnarray*}
&&\p \Biggl( \Biggl\llVert \frac{1}{n}\sum_{i=0}^{n-1}
g(X_i) - \bfH \Biggr\rrVert _{\mathrm{HS}} \ge C\sqrt {
\frac{m}{n}}\tau \llVert g\rrVert _\infty+ t \Biggr)
\\
&&\quad \le 2\exp \biggl(-cn\min \biggl(\frac{t^2}{m\llVert g\rrVert _\infty^2 \tau
^2},\frac{t}{m\llVert g\rrVert _\infty\tau\log\tau}
\biggr) \biggr).
\end{eqnarray*}
It is easy to see that by adjusting the value of the absolute constant
this is equivalent to
%
\begin{eqnarray}
\label{eq:final_estimate_sum} \p \Biggl( \Biggl\llVert \frac{1}{n}\sum
_{i=0}^{n-1} g(X_i) - \bfH \Biggr\rrVert
_{\mathrm{HS}} \ge t \Biggr) \le 2\exp \biggl(-cn\min \biggl(\frac{t^2}{m\llVert g\rrVert _\infty^2 \tau
^2},
\frac{t}{m\llVert g\rrVert _\infty\tau\log\tau} \biggr) \biggr).
\end{eqnarray}
Since by \eqref{eq:L_infty_bound} $\llVert g\rrVert _\infty= \sup_{x \in\calX}
|h(x,x)|$, to finish the proof of Theorem~\ref{thm:exponential_ineq}
it is enough to combine the above inequality with Theorem~\ref{thm:HW}.
\end{pf*}

\begin{rem} Let us mention that a Markov chain is geometrically ergodic
iff it satisfies the following drift condition (see Theorem~16.0.1. in
\cite{MT}). There exists $\lambda\in(0,1)$, $b \in\R_+$ and
$V\dvtx \mathcal{X} \to[1,\infty)$ such that for some set $C$,
satisfying \eqref{eq:small_set} and $\pi(C) > 0$,
\begin{eqnarray*}
P^m V - V \le-\lambda V + b\Ind{C}
\end{eqnarray*}
and $K := \sup_{x\in C} V(x) < \infty$. Finding appropriate drift
functions is in fact the most common way of proving geometric ergodicity.

It turns out that one can bound the quantity $\tau$ appearing in the
estimate \eqref{eq:final_estimate_sum} in terms of the parameters of
the drift condition and \eqref{eq:small_set}. Such an estimate follows
directly from Propositions 6, 7 from \cite{AdBed} (obtained with help
of previous important estimates from \cite{Bax}). Namely for a chain
started from a point $x$, we have
\begin{eqnarray*}
\tau&\le& 2\log \biggl(\frac{\log(\afrac{6}{2-\delta})}{\log(\afrac
{6}{2-\delta})} \biggr)
\\
&&{}\times\max \biggl(\frac{\log(V(x)\Ind{C^c}(x)+(b(1-\lambda
)^{-1}+K)\Ind{C}(x))}{\log2},\frac{\log(b(1-\lambda)^{-1}+K)}{\log
2},1 \biggr)
\\
&&{}\times \frac{1}{\log\sklafrac{1}{1-\lambda}},
\end{eqnarray*}
where $\delta$ is the parameter from \eqref{eq:small_set}.
\end{rem}

%
\begin{rem}\label{re:starting_measure} It is clearly seen from the
proof of Theorem~\ref{thm:exponential_ineq} that the chain does not
have to be started from a point. It is sufficient to assume that the
stopping time $T_0$ is exponentially integrable under the starting
measure $\mu$. This will be the case, for example, if the function $V$
in the drift conditions is $\mu$-integrable (as follows by Proposition~4.1. (ii) in \cite{Bax}).
\end{rem}

%
\begin{rem}
We also note that the absolute constant $c$ in \eqref
{eq:final_estimate_sum} can be given explicitly, since
Lemma~\ref{le:Bernstein_psi_one} with explicit constants is known \cite{VershyninFrame}, the
constant from Lemma~3 in \cite{AdMarkovTail} can be easily read from
the proof and the L\'evy type inequality by Montgomery--Smith is also
given with explicit constants \cite{MontgomerySmithLevy}. We do not
pursue this direction here. See \cite{AdBed} for related inequalities
for additive functionals of Markov chains with explicit constants.
\end{rem}

%
\section{Discussion of optimality. Counterexamples}
\label{sec:examples}
We would like to conclude with an example of a square integrable kernel
$h$ and a uniformly ergodic Markov chain for which the conclusion of
Theorem~\ref{thm:LLN_spectrum} fails and the empirical counterpart of
the spectrum almost surely is not convergent to the spectrum of $\bfH
$. The example uses directly the construction of \cite{ABDGHW}, where
a counterexample to the law of large numbers for $U$-statistics was
given. We adapt it to our setting and provide the details for the sake
of completeness.

Let thus $\varepsilon_0,\varepsilon_1,\ldots$ be i.i.d. random
variables with distribution $\p(\varepsilon_i = 1) = \p(\varepsilon
_i = 0) = 1/2$ and $Y_0,Y_1,\ldots$ -- i.i.d. random variables
distributed uniformly on the interval $(0,1)$, independent of the
sequence $(\varepsilon_i)$. Define $X_0 = x$,
\begin{eqnarray*}
X_{n+1} = \lleft\{ %
\begin{array} {l@{\qquad}l}
X_n& \mbox{if } \varepsilon_n = 0,
\\
Y_{n+1}& \mbox{if } \varepsilon_n = 1. \end{array}
\rright.
\end{eqnarray*}
It is easy to see that \eqref{eq:small_set} is satisfied with $m=1$,
$\delta= 1/2$, $C = (0,1)$ and $\nu$ being the Lebesgue measure on
$(0,1)$, thus (see \cite{MT}) the chain is uniformly ergodic (i.e., it
is geometrically ergodic and the function $M(x)$ in \eqref
{eq:geom_erg} is bounded by a constant independent of $x$). The unique
stationary measure for the chain, $\pi$ is in this case the Lebesgue
measure. Consider now a function $h\dvtx (0,1)^2 \to\R$ given by
$h(x,y) = 0$ if $x\neq y$ and $h(x,x) = 1/x^3$.

Of course $h = 0$ $\pi\otimes\pi$-a.s. and so $\bfH= 0$. Let now
$i_0<i_1<i_2<\cdots$ be defined as $i_0 = 0$, $i_{n+1} = \min\{ i >
i_n\dvt \varepsilon_i = 0,\varepsilon_{i+1}=1\}$. Then for $k > 0$,
$X_{i_k} = X_{i_{k}+1}$, moreover conditionally on $(\varepsilon
_i)_{i\ge0}$, $X_{i_k}$ are i.i.d., distributed according to $\pi$.
Since the absolute value of the largest eigenvalue of a matrix is not
smaller than the absolute value of its maximal entry, both $\bfH_n$
and $\tbfH_n$ have at least one eigenvalue, which in absolute value exceeds
$n^{-1}\max_{0\le i\le n-2} h(X_i,X_{i+1})$. Moreover, by the law of
large numbers $i_n/n \to4$ a.s., so using the conditional independence
of $X_i$, the Borel--Cantelli lemma and the fact that $\p(X_{i_k} \le
t) = t$ for $t \in(0,1)$, we get
\begin{eqnarray*}
\limsup_{n\to\infty} \max\lambda(\bfH_n), \limsup
_{n\to\infty} \max\lambda(\tbfH_n) &\ge&\limsup
_{k\to\infty} \frac
{1}{i_k+2}h(X_{i_k},X_{i_k+1})
\\
&=& \limsup_{k\to\infty} \frac{1}{i_k+2}\frac{1}{X_{i_k}^3} =
\infty \qquad \mbox{a.s.}
\end{eqnarray*}

This shows that the law of large numbers for spectra fails in this case.


\section*{Acknowledgements}
Research partially supported by NCN Grant N N201 608740.
We would like to thank the anonymous referee for all the remarks
concerning the submitted version of the manuscript.



\printhistory
\end{document}